\documentclass{amsart}

\newtheorem{theorem}{Theorem}
\newtheorem*{theorem*}{Theorem}
\newtheorem{lemma}{Lemma}
\newtheorem*{lemma*}{Lemma}

\newtheorem*{corollary*}{Corollary}

\newtheorem*{proposition*}{Proposition}

\newtheorem*{conjecture*}{Conjecture}
{\theoremstyle{definition} \newtheorem{definition}{Definition}
	\newtheorem*{definition*}{Definition}
	\newtheorem{example}{Example}
	\newtheorem*{example*}{Example}
	
	\newtheorem*{remark*}{Remark}
	
	\newtheorem*{note*}{Note}
}

\usepackage[
pdftex,
colorlinks,%
linkcolor=blue,citecolor=red,urlcolor=blue,
hyperindex,%
plainpages=false,%
bookmarksopen,%
bookmarksnumbered%
]{hyperref}

\begin{document}

	\title[Integration of the nonlinear Schr\"{o}dinger equation ...]{Integration of the nonlinear Schr\"{o}dinger equation with a self-consistent source and nonzero boundary conditions}
	\author{Anvar Reyimberganov}
	
	\address{Urgench State University, Uzbekistan}
	\email{anvar@urdu.uz}
	\urladdr{http://www.urdu.uz}

\maketitle

	\begin{abstract}
		This paper is devoted to the study of the defocusing nonlinear Schr\"{o}dinger equation with a self-consistent source and nonzero boundary conditions by the method of the inverse scattering problem. In cases where the source consists of a combination of eigenfunctions of the corresponding spectral problem for the Zakharov-Shabat system, the complete integrability of the nonlinear Schr\"{o}dinger equation is investigated. Namely, the evolutions of the scattering data of the self-adjoint Zakharov-Shabat system, whose potential is a solution of the defocusing nonlinear Schr\"{o}dinger equation with a self-consistent source, are obtained.
	\end{abstract}
\section{Introduction}

Nonlinear Schr\"{o}dinger (NLS) equation
\begin{equation} \label{eq1} 
	iu_{t} -2\chi \left|u\right|^{2} u+u_{xx} =0 
\end{equation} 
with various boundary conditions models a wide class of nonlinear phenomena in physics. In the work \cite{Zakharov1971}, V. Zakharov and A. Shabat showed that NLS equation can be applied in the study of optical self-focusing and splitting of optical beams. The sign of the coupling constant corresponds to the attraction ($\chi <0$) and repulsion ($\chi >0$) of particles. In the case of attraction, the problem of a finite number of particles and their bound states has a physical meaning. In the classical limit, this is modeled by rapidly decreasing boundary conditions. In the case of repulsion, the problem corresponding to a gas of particles with a finite density is of interest. 

It is well known that inverse scattering method for integration of the NLS equation is based on the so-called Zakharov-Shabat and Ablowitz-Kaup-Newell-Segur scatte\-ring problem (see: \cite{Ablowitz1974, Zakharov1971, Faddeev1987}). 

In 1971, V. Zakharov and A. Shabat \cite{Zakharov1971} showed that the NLS equation can be solved by means of the inverse scattering transform (IST) technique. The IST as a method to solve the initial-value problem for the NLS equation has been extensively studied in the literature, both in the focusing ($\chi =-1$) and in the defocusing ($\chi =1$) dispersion regimes. The IST for the defocusing NLS equation with nonzero boundary conditions was first studied in 1973 by V.~Zakharov and A. Shabat \cite{Zakharov1972} and a detailed study can be found in the monograph \cite{Faddeev1987}. 

In the work \cite{Demontis2013} a rigorous theory of the IST for the defocusing nonlinear Schr\"{o}dinger equation with nonvanishing boundary values $u_{\pm } =u_{0} e^{i\theta _{\pm } } $ as $x\to \pm \infty $ is presented. The IST theory for the defocusing NLS equation with nonzero boundary conditions was studied by B. Gino, F. Emily and B. Prinari \cite{Biondini2016} and the focusing case has been studied by G. Biondini, G. Kova\v{c}i\'{c} \cite{Biondini2014}, F. Demontis, B. Prinari, C. van der Mee, F. Vitale \cite{Demontis2014}. 

V.K. Melnikov \cite{Melnikov1991, Melnikov1992} showed that the NLS equation remains its integrability by the inverse scattering method,  if a source is added to them in the form of a combination of eigenfunctions of the corresponding spectral problem. Namely, the term ``self-consistent source'' was introduced in the works of V.K. Melnikov.

The NLS equation with the self-consistent sources in various classes of functions were considered by A.B. Khasanov \cite{Khasanov2009}, I.D. Rakhimov \cite{Reyimberganov2021}, A.B. Yakhshimuratov \cite{Yakhshimuratov2011}. 

In the matrix case, the inverse scattering theory for the matrix Zakharov-Shabat system  was investigated by P. Barbara, F. Demontis and C. Van der Mee \cite{Prinari2018, Demontis2014Acta, Demontis2022} and applied for the integration of the matrix NLS equation. In \cite{Urazboev2021}, the matrix NLS equation with the self-consistent source was considered in the class of rapidly decreasing matrix functions.

\section{Formulation of the problem}

We consider the following system of equations 
\begin{gather}
	iu_{t} -2u\left|u\right|^{2} +u_{xx} =-2i\sum _{n=1}^{N}(f_{1,n}^{*} g_{2,n}^{*} +f_{2,n} g_{1,n} )   \label{eq11} 
	\\
	\frac{\partial f_{1,n} }{\partial x} -u^{*} f_{2,n} +i\xi _{n} f_{1,n} =\frac{\partial f_{2,n} }{\partial x} -uf_{1,n} -i\xi _{n} f_{2,n} =0,\, \, \, n=1,2,...,N, \label{eq12}
	\\
	\frac{\partial g_{1,n} }{\partial x} -ug_{2,n} -i\xi _{n} g_{1,n} =\frac{\partial g_{2,n} }{\partial x} -u^{*} g_{1,n} +i\xi _{n} g_{2,n} =0,\, \, \, n=1,2,...,N \label{eq13}  
\end{gather} 
under the initial condition
\begin{equation} \label{eq14} 
	u(x,0)=u_{0} (x),\quad x\in R.  
\end{equation} 
Here the initial function $u_{0} (x)$, $ x\in R$ satisfies the following properties:

1) 
\begin{equation} \label{eq15} 
	\int _{-\infty }^{0}(1-x)\left|u_{0} (x)-\rho e^{i\alpha _{-} } \right|dx\,  +\int _{0}^{\infty }(1+x)\left|u_{0} (x)-\rho e^{i\alpha _{+} } \right|dx<\infty  ,\,  
\end{equation} 
where $\rho >0$ and $0\le \alpha _{\pm } <2\pi $ are arbitrary constants.

2) The system of equations \eqref{eq12} with coefficient $u_{0} (x)$ possesses exactly $N$ eigenvalues $\xi _{1} (0)$, $\xi _{2} (0)$, ... , $\xi _{N} (0)$. 

We assume that the solution $u(x,t)$ of the equation \eqref{eq11} is sufficiently smooth and tends to its limits sufficiently  rapidly as $x\to \pm \infty $, i.e., for all $t\ge 0$ satisfies  the condition
\begin{align} 
	\int _{-\infty }^{0}(1-x)\left|u(x,t)-\rho e^{i\alpha _{-} -2i\rho ^{2} t} \right|dx+&\int _{0}^{\infty }(1+x)\left|u(x,t)-\rho e^{i\alpha _{+} -2i\rho ^{2} t} \right|dx  \nonumber\\
	+&\int _{-\infty }^{\infty }\sum _{k=1}^{2}\left|\frac{\partial ^{k} u(x,t)}{\partial x^{k} } \right| dx<\infty.  \label{eq16}
\end{align}
It follows from this condition that the left-hand side of equation \eqref{eq11} for all $t\ge 0$ tends to zero as $x\to \pm \infty $. Taking this into account, the solutions $\, F_{n} (x,t)=(f_{1,n} ,f_{2,n} )^{T}  $ and $G_{n} (x,t)=(g_{1,n} ,g_{2,n} )^{T} $ of equations \eqref{eq12} and \eqref{eq13}, respectively, are to be chosen so that the expression in the right-hand side of equation \eqref{eq11} for all $t\ge 0$ should tend to zero rapidly enough as $x\to \pm \infty $. This can be done in the following two ways:

(A) Let the functions $F_{n} (x,t)$ and $G_{n} (x,t)$ be the eigenfunctions of equations \eqref{eq12} and \eqref{eq13} respectively, corresponding to the eigenvalues $\xi _{n} (t)$. In this case, the functions $F_{n} (x,t)$ and $G_{n} (x,t)$ belong to the $L^{2} (R)$ for all $ t\ge 0$. We assume that
\begin{equation} \label{eq17} 
	\int _{-\infty }^{\infty }G_{n}^{T} (s,t)F_{n} (s,t)ds=A_{n} (t) ,\quad t\ge 0,\quad n=1,2,...,N.  
\end{equation} 
Where $A_{n} (t)$ are given and the continuous functions of $t$.

(B) Let the function $F_{n} (x,t)$ be the eigenfunctions of equations \eqref{eq12} corresponding to the eigenvalues $\xi _{n} (t)$ and let the function $G_{n} (x,t)$ be unbounded solution of the equation \eqref{eq13}, satisfying equalities 
\begin{equation} \label{eq18} 
	f_{1,n} g_{1,n} -f_{2,n} g_{2,n} =B_{n} (t),\quad t\ge 0,\quad n=1,2,...,N.  
\end{equation} 
Where the functions $B_{n} (t)$ are given continuous real valued functions of $t$.

Our main goal is to obtain representations for solutions $u(x,t)$, $F_{n} (x,t)$, $G_{n} (x,t)$, $n=1,2,...,N$ of problem \eqref{eq11}-\eqref{eq16}, within the framework of the inverse scattering method.

\section{Preliminaries}

Let a function $u(x,t)$ belong to the class of functions \eqref{eq16}. In this section, we give well known \cite{Faddeev1987}, necessary information concerning the theory of direct and inverse scattering problems for the operator 
\[L(t)=i\left(\begin{array}{cc} {\frac{\partial }{\partial x} } & {-u^{*} (x,t)} \\ {u(x,t)} & {-\frac{\partial }{\partial x} } \end{array}\right), \quad t\ge 0.\] 

Consider the equation
\begin{equation} \label{ZEqnNum935108} 
	(L(t)-\xi I)f=0 
\end{equation} 
with respect to an unknown $2\times 2$ square matrix function $f(x,\xi ,t)$. Here $\xi $\textit{ }is a complex spectral parameter. 

We introduce a new spectral parameter $p$ as follows $p=\sqrt{\xi ^{2} -\rho ^{2} } $. The variable $p$ is then thought of as belonging to a Riemann surface $\Gamma $ consisting of a sheet $\Gamma _{+} $ and a sheet $\Gamma _{-}$ which both coincide with the complex plane cut along the semi lines $\Sigma =(-\infty ,-\rho ]\cup [\rho ,\infty )$ with its edges glued in such a way that $p(\xi )$ is continuous through the cut. The variable $p$ is thought of as belonging to the complex plane consisting of the upper half complex plane $\Gamma _{+} $ and the lower half complex plane $\Gamma _{-} $ glued together along the whole real line. For all $\xi \in \Sigma $, the branch of the square root is fixed by the condition $sign\, p(\xi )=sign\, \xi $.

For $\xi \in \Sigma $, we define matrix Jost solutions $f^{-} (x,\xi ,t)$ and $f^{+} (x,\xi ,t)$ from the right and the left, respectively as those square matrix solutions to \eqref{ZEqnNum935108} satisfying asymptotics
\begin{equation}\label{ZEqnNum770859}
	f^{\pm} \sim E^{\pm}(x,\xi,t)   \quad \text{as} \quad x \to \pm \infty
\end{equation}
where
\[E^{\pm } (x,\xi ,t)=\left(\begin{array}{cc} {1} & {-\frac{i(\xi -p)}{\rho } e^{-i\alpha _{\pm } +2i\rho ^{2} t} } \\ {\frac{i(\xi -p)}{\rho } e^{i\alpha _{\pm } -2i\rho ^{2} t} } & {1} \end{array}\right)e^{-ip\sigma _{3} x} .\] 
Here and everywhere below we will use the standard Pauli matrices
\[\sigma _{1} =\left(\begin{array}{cc} {0} & {1} \\ {1} & {0} \end{array}\right), \quad \sigma _{2} =\left(\begin{array}{cc} {0} & {-i} \\ {i} & {0} \end{array}\right), \quad \sigma _{3} =\left(\begin{array}{cc} {1} & {0} \\ {0} & {-1} \end{array}\right).\] 
If a function $u(x,t)$ belongs to the class of functions \eqref{eq16}, then such a solution to equations \eqref{ZEqnNum935108} exists and is unique.

It can be shown that
\begin{equation} \label{ZEqnNum112754} 
	\frac{d}{dx} \det f_{}^{\pm } (x,\xi ,t)=0.  
\end{equation} 
From \eqref{ZEqnNum112754} and \eqref{ZEqnNum770859} it follows that 
\begin{equation} \label{ZEqnNum736324} 
	\det f_{}^{\pm } (x,\xi ,t)=\frac{2p(\xi -p)}{\rho ^{2} } .  
\end{equation} 
The system \eqref{ZEqnNum935108} is invariant with respect to the involution
\begin{equation} \label{ZEqnNum603247} 
	\bar{f}^{\pm } (x,\xi ,t)=\sigma _{1} f^{\pm } (x,\xi ,t)\sigma _{1},\quad \xi \in \Sigma .  
\end{equation} 
We now call the columns of
\[
f^{+} (x,\xi ,t)=\left(\bar{\psi }(x,\xi ,t)\, \, \, \, \psi (x,\xi ,t)\right),
\,\,\,
f^{-} (x,\xi ,t)=\left(\varphi (x,\xi ,t)\, \, \, \, \bar{\varphi }(x,\xi ,t)\right)
\] 
the Jost solutions from the right and the left, respectively. For the Jost solutions we get the following asymptotic estimates
\begin{equation} \label{ZEqnNum260418} 
	\psi (x,\xi ,t)\sim \left(\begin{array}{c} {-\frac{i(\xi -p)}{\rho } e^{-i\alpha _{+} +2i\rho ^{2} t} } \\ {1} \end{array}\right)e^{ipx}, \quad x\to \infty,  
\end{equation} 
\[
\bar{\psi }(x,\xi ,t)\sim \left(\begin{array}{c} {1} \\ {\frac{i(\xi -p)}{\rho } e^{i\alpha _{+} -2i\rho ^{2} t} } \end{array}\right)e^{-ipx} ,\quad  x\to \infty ,  
\]
\begin{equation} \label{ZEqnNum830052} 
	\varphi (x,\xi ,t)\sim \left(\begin{array}{c} {1} \\ {\frac{i(\xi -p)}{\rho } e^{i\alpha _{-} -2i\rho ^{2} t} } \end{array}\right)e^{-ipx} ,\quad x\to -\infty,  
\end{equation} 
\[	
\bar{\varphi }(x,\xi ,t)\sim \left(\begin{array}{c} {-\frac{i(\xi -p)}{\rho } e^{-i\alpha _{-} +2i\rho ^{2} t} } \\ {1} \end{array}\right)e^{ipx} ,\quad  x\to -\infty,  
\]

Since $f^{-} (x,\xi ,t)$ and $f^{+} (x,\xi ,t)$ are square matrix solutions of the homogeneous first order equation \eqref{ZEqnNum935108}, we necessarily have for $\xi \in \Sigma $
\begin{equation} \label{ZEqnNum733217} 
	f^{-} (x,\xi ,t)=f^{+} (x,\xi ,t)S(\xi ,t) 
\end{equation} 
where $S(\xi ,t)$ is the transition coefficient matrix. 

From the involution property \eqref{ZEqnNum603247} for $\xi \in \Sigma $, it follows that
\[\bar{S}(\xi ,t)=\sigma _{1} S(\xi ,t)\sigma _{1} .\] 
Hence, we have
\begin{equation} \label{ZEqnNum320925} 
	S(\xi ,t)=\left(\begin{array}{cc} {a(\xi ,t)} & {\bar{b}(\xi ,t)} \\ {b(\xi ,t)} & {\bar{a}(\xi ,t)} \end{array}\right).  
\end{equation} 
Coefficients $a(\xi ,t)$ and $b(\xi ,t)$ are called scattering coefficients. From relations \eqref{ZEqnNum736324} and \eqref{ZEqnNum320925} we obtain
\[
a(\xi ,t)\bar{a}(\xi ,t)-b(\xi ,t)\bar{b}(\xi ,t)=1.  
\]
By using \eqref{ZEqnNum733217} we can represent the scattering coefficients as
\begin{equation} \label{ZEqnNum762563} 
	a(\xi ,t)=\frac{\rho ^{2} }{2p(\xi -p)} \det (\varphi (x,\xi ,t),\, \psi (x,\xi ,t)) 
\end{equation} 
and
\[b(\xi ,t)=\frac{\rho ^{2} }{2p(\xi -p)} \det (\bar{\psi }(x,\xi ,t),\varphi (x,\xi ,t)).\] 

If a function $u(x,t)$ belongs to the class of functions \eqref{eq16}, then for each $x\in R$ the Jost solutions $\psi (x,\xi ,t)e^{-ipx} $ and $\varphi (x,\xi ,t)e^{ipx} $ are analytic for $\xi \in \Gamma _{+} $ excluding branch points $\xi =\pm \rho $, there are asymptotes for $\left|\xi \right|\to \infty $
\begin{equation} \label{ZEqnNum910868} 
	\varphi (x,\xi ,t)e^{ipx} =\left(\begin{array}{c} {1} \\ {\frac{i(\xi -p)}{\rho } e^{i\alpha_{-} -2ip^{2} t} } \end{array}\right)+O\left(\frac{\left|1+\xi -p\right|}{\left|\xi \right|} \right),  
\end{equation} 
\begin{equation} \label{ZEqnNum318020} 
	\psi (x,\xi ,t)e^{-ipx} =\left(\begin{array}{c} {-\frac{i(\xi -p)}{\rho } e^{-i\alpha_{+} +2ip^{2} t} } \\ {1} \end{array}\right)+O\left(\frac{\left|1+\xi -p\right|}{\left|\xi \right|} \right).  
\end{equation} 

It follows from the analyticity properties of the Jost solutions and equality \eqref{ZEqnNum762563} that the function $a(\xi ,t)$ can be analytically continued to the sheet $\Gamma _{+} $ excluding branch points $\xi =\pm \rho $. 

From \eqref{ZEqnNum910868} and \eqref{ZEqnNum318020} we obtain that for $\left|\xi \right|\to \infty $, the function $a(\xi ,t)$ has the asymptotics
\begin{equation} \label{ZEqnNum447933} 
	a(\xi ,t)=1+O\left(\frac{1}{\left|\xi \right|} \right)\, \, \, {\rm as}\, \, Im\, \xi >0 
\end{equation} 
and
\begin{equation} \label{ZEqnNum687898} 
	a(\xi ,t)=e^{-i\theta } +O\left(\frac{1}{\left|\xi \right|} \right)\, \, {\rm as}\, \, Im\, \xi <0 
\end{equation} 
where we recall $\theta =\alpha _{+} -\alpha _{-} $. 

Similarly, the function $\bar{a}(\xi ,t)$ can be analytically continued to the sheet $\Gamma _{-} $, excluding branch points $\xi =\pm \rho $.

It follows from the analyticity with respect to the function $a(\xi ,t)$ on $\Gamma _{+} $ and from the asymptotics \eqref{ZEqnNum447933}, \eqref{ZEqnNum687898} that the function $a(\xi ,t)$ can have only a finite number of zeros on the sheet $\Gamma _{+} $. These zeros will be denoted by $\xi _{1} ,\xi _{2} ,...,\xi _{N} $. In [4] it is shown that all zeros are simple and all belong to the $(-\rho ,\rho )$.

It is seen from  representation \eqref{ZEqnNum762563}  that $\xi =\xi _{n} $ the functions  $\varphi (x,\xi ,t)$ and $\psi (x,\xi ,t)$ are proportional to each other
\begin{equation} \label{ZEqnNum197225} 
	\varphi _{n} (x,t)=c_{n} (t)\psi _{n} (x,t),\, \, \, \, \, \, \, \bar{\varphi }_{n} (x,t)=c_{n}^{*} (t)\bar{\psi }_{n} (x,t),\, \, \, \, \, n=1,2,...,N.  
\end{equation} 
Where $\varphi _{n} (x,t)=\varphi (x,\xi _{n} ,t),\, \, \, \psi _{n} (x,t)=\psi (x,\xi _{n} ,t)$. 

The zeros of $a(\xi,t)$ correspond to the eigenvalues of the equation \eqref{ZEqnNum935108}.  The equation \eqref{ZEqnNum935108} is self-adjoint, so its eigenvalues and thus the zeros of the function $a(\xi ,t)$ are real.

Note, the vector functions
\begin{equation} \label{ZEqnNum510514} 
	h_{n} (x,t)=\frac{\frac{\partial }{\partial \xi } \left. \left(\varphi (x,\xi ,t)-c_{n} \psi (x,\xi ,t)\right)\right|_{\xi =\xi _{n} } }{\dot{a}(\xi _{n} ,t)} ,\, \, n=1,2,...,N 
\end{equation} 
are a solution to the equations $(L(t)-\xi _{n} I)h_{n} =0$. Where $\dot{a}(\xi _{n} ,t)=\frac{\partial }{\partial \xi } \left. a(\xi ,t)\right|_{\xi =\xi _{n} } $.

From equality \eqref{ZEqnNum510514} it follows that
\[	h_{n} (x,t)\sim -c_{n} (t)\left(\begin{array}{c} {-\frac{i(\xi _{n} -p_{n} )}{\rho } e^{-i\alpha _{-} +2ip^{2} t} } \\ {1} \end{array}\right)e^{ip_{n} x}  \,\,as \,\,x\to -\infty , 
\]
\begin{gather}
	h_{n} (x,t)\sim \left(\begin{array}{c} {1} \\ {\frac{i(\xi _{n} -p_{n} )}{\rho } e^{i\alpha _{+} -2ip^{2} t} } \end{array}\right)e^{-ip_{n} x}  \,\,as \,\,x\to \infty, \label{ZEqnNum922580}
\end{gather}
where $p_{n} =i\sqrt{\rho ^{2} -\xi _{n}^{2} } $. In particular, we have
\begin{equation} \label{ZEqnNum819680} 
	\det (\varphi _{n} ,h_{n} )=-\frac{2p_{n} (\xi _{n} -p_{n} )}{\rho ^{2} } c_{n} ,\, \, n=1,2,...,N.  
\end{equation} 

\begin{definition}\label{def1}
	The set $\{ a(\xi ,t),\, \, b(\xi ,t),\, \, \xi _{n} (t),\, \, c_{n} (t),\, \, n=1,2,...,N\} $ is called the scattering data for equation \eqref{ZEqnNum935108}. The direct scattering problem is to find the scattering data via the given potentials $u(x,t)$ and the inverse scattering problem is to find the potentials $u(x,t)$ of the equation \eqref{ZEqnNum935108} via the given scattering data.
\end{definition}

Before we proceed further solving the inverse problem, it is convenient to introduce a uniformization variable $z$ (see \cite{Demontis2013, Faddeev1987}) defined by the conformal mapping: $z=z(\xi )=\xi +p(\xi )$. Inverse mapping given by
\[\xi =\frac{1}{2} \left(z+\frac{\rho ^{2} }{z} \right), \,\,\,p=z-\xi =\frac{1}{2} \left(z-\frac{\rho ^{2} }{z} \right).\] 
With this mapping the sheets $\Gamma _{+} $ and $\Gamma _{-} $ of the Riemann surface $\Gamma $ are, respectively, mapped onto the upper and lower complex half-planes $Im\, z>0$ and $Im\, z<0$ of the complex $z$-plane. The cut $\Sigma $ on the Riemann surface is mapped onto the real $z$ axis. The segments $[-\rho ,\rho ]$ on $\Gamma _{+} $ and $\Gamma _{-} $ are mapped onto the upper and lower semicircles of radius $\rho $ and center at the origin of the $z$-plane. The neighborhood of the point $\xi =\infty $ on $\Gamma _{\pm } $ with the condition $\pm Im\xi >0$ is mapped into the neighborhood of the point $z=\infty $, and the neighborhood of the point $\xi =\infty $ on $\Gamma _{\pm } $ with the condition $\pm Im\xi <0$ is mapped into the neighborhood of the point $z=0$. 

In terms of variable $z$, relation \eqref{ZEqnNum733217} can be written when $Im z=0$ following form
\begin{equation} \label{39)} 
	f^{-} (x,z,t)=f^{+} (x,z,t)S(z,t) 
\end{equation} 
here $f^{\pm } (x,z,t)=f^{\pm } (x,\xi (z),t)$, $S(z,t)=S(\xi (z),t)$ and one can obtain the symmetries of the scattering coefficients:
\begin{equation} \label{ZEqnNum241555} 
	a(z,t)=\bar{a}\left(\frac{\rho ^{2} }{z} ,t\right),\, \, \, Im z\ge 0,  
\end{equation} 
\[
b(z,t)=-\bar{b}\left(\frac{\rho ^{2} }{z} ,t\right),\, \, \, Im z=0.  
\] 

Equality \eqref{ZEqnNum241555}, together with the self-adjointness of the equation \eqref{ZEqnNum935108}, ensure that the scattering coefficient $a(z,t)$ ($\bar{a}(z,t)$) can only have zeros at $z_{n} =\xi _{n} +iv_{n} $  ($\bar{z}_{n} =\xi _{n} -iv_{n} $), with $-\rho <\xi _{n} <\rho $ and $v_{n} =\sqrt{\rho ^{2} -\xi _{n}^{2} } >0$.

Taking into account the analyticity properties of $a(z,t)$ in the upper half plane $Im z>0 $  we can obtain the following representation
\[
a(z,t)=\prod_{n=1}^{N}\frac{z-z_n}{z-z_{n}^{*}} \exp\left[-\frac{1}{2 \pi i} \int_{-\infty }^{\infty }\frac{\log(1-|r(\zeta,t)|^2)}{\zeta-z}d\zeta \right].
\]
Where $r(z,t)\equiv \frac{b(z,t)}{a(z,t)}$ is called reflection coefficient. 
According to \eqref{ZEqnNum687898}, for $z \to 0$ we obtain that
\[
e^{-i\theta }=\prod_{n=1}^{N}\frac{z_n}{z_{n}^{*}} \exp\left[-\frac{1}{2 \pi i} \int_{-\infty }^{\infty }\frac{\log(1-|r(\zeta,t)|^2)}{\zeta}d\zeta \right].
\]

If $\varphi _{n} =\left(\begin{array}{l} {\varphi _{1,n} } \\ {\varphi _{2,n} } \end{array}\right)$ is an eigenfunction of the equation \eqref{ZEqnNum935108} corresponding to $z_{n} $, then we define $\bar{\varphi }_{n} =\left(\begin{array}{l} {\varphi _{2,n}^{*} } \\ {\varphi _{1,n}^{*} } \end{array}\right)$ to be the eigenfunction of the equation \eqref{ZEqnNum935108} corresponding to $\bar{z}_{n} $.

It is well known that the inverse scattering theory of  \eqref{ZEqnNum935108} can be formulated in terms of the Gelfand-Levitan-Marchenko equations. The Jost solution $\psi (x,z,t)$ of the equation \eqref{ZEqnNum935108} can be represented in the following form
\begin{equation} \label{ZEqnNum820131} 
	{\psi (x,z ,t)=\left(\begin{array}{c} {-\frac{i\rho }{z} e^{-i\alpha _{+} +2i\rho ^{2} t} } \\ {1} \end{array}\right)e(x,z)+\int _{x}^{\infty }\mathcal{K} (x,y,t)\left(\begin{array}{c} {-\frac{i\rho }{z} e^{-i\alpha _{+} +2i\rho ^{2} t} } \\ {1} \end{array}\right)e(y,z)dy ,} 
\end{equation} 
here $e(x,z)=e^{\frac{i}{2} \left(z-\frac{\rho ^{2} }{z} \right)x} $ and $\mathcal{K}(x,y)$ is a 2 $\mathrm{\times}$ 2 matrix function which has to satisfy the following Gelfand-Levitan-Marchenko equation:
\[
\mathcal{K}(x,y,t)+\mathcal{F}(x+y,t)+\int _{x}^{\infty }\mathcal{K}(x,s,t)\mathcal{F}(s+y,t)ds =0,\, \, \, y\ge x,  
\]
where $\mathcal{K}(x,y,t)$ and $\mathcal{F}(x+y,t)$ are defined as 
\[	\mathcal{K}(x,y,t)=\left(\begin{array}{cc} {K_{11}^{} (x,y,t)} & {K_{12}^{} (x,y,t)} \\ {K_{21} (x,y,t)} & {K_{22} (x,y,t)} \end{array}\right), \,\, \mathcal{F}(x+y,t)=\left(\begin{array}{cc} {F_{1} (x+y,t)} & {F_{2}^{*} (x+y,t)} \\ {F_{2} (x+y,t)} & {F_{1} (x+y,t)} \end{array}\right)\] 
with
\[F_{1} (x,t)=\frac{\rho e^{i\alpha_{+} -2i\rho ^{2} t} }{4\pi i} \int _{-\infty }^{\infty }\frac{b(z,t)}{za(z,t)} \cdot e(x,z)dz -\frac{1}{2} \sum _{n=1}^{N}\frac{c_{n} (t)\rho e^{-i\alpha_{+} +2i\rho ^{2} t} }{\dot{a}(z_{n} ,t)z_{n} }  \cdot e(x,z_{n} ),\] 
\[F_{2} (x,t)=\frac{1}{4\pi } \int _{-\infty }^{\infty }\frac{b(z,t)}{a(z,t)} \cdot e(x,z)dz -\frac{1}{2} \sum _{n=1}^{N}\frac{ic_{n} (t)}{\dot{a}(z_{n} ,t)}  \cdot e(x,z_{n} ).\] 
In representations \eqref{ZEqnNum820131}, the component $K_{21}^{} (x,x,t)$ of the matrix $\mathcal{K}(x,y,t)$ have relations with the potential 
\[
2K_{21}(x,x,t)=\rho e^{i\alpha _{+} -2i\rho ^{2} t} -u(x,t).  
\]
In the work \cite{Faddeev1987}, it was proven the uniquely determining of the potential $u(x,t)$ by the scattering data.

\section{Time evolution}

The use of the inverse scattering method for integration of the problem \eqref{eq11}-\eqref{eq16} is based on the following. Let the function $u(x,t)$ be a solution of equation \eqref{eq11}, from the class of functions \eqref{eq15}. Consider equation \eqref{ZEqnNum935108} with a potential $u(x,t)$ and find the evolution from the scattering data.

Assuming that
\begin{equation} \label{ZEqnNum767561} 
	F_{N+n} =\left(\begin{array}{l} {f_{2,n}^{*} } \\ {f_{1,n}^{*} } \end{array}\right),\, \, \, G_{N+n} =\left(\begin{array}{l} {g_{2,n}^{*} } \\ {g_{1,n}^{*} } \end{array}\right),\, \, \, \, n=1,2,...,N,  
\end{equation} 
equation \eqref{eq11} can be represented as an equality of operators in the class of smooth functions $f(x,\xi ,t)$ satisfying the equation \eqref{ZEqnNum935108}: 
\[ 
\frac{\partial L}{\partial t} +[L,A]=i\sum _{n=1}^{2N}[\sigma _{3} ,\, F_{n}^{} G_{n}^{T} ] .  
\]
Where $[L,A]=LA-AL$ and 
\[A=\left(\begin{array}{cc} {i\left|u\right|^{2} +2i\xi ^{2} } & {-iu_{x}^{*} -2\xi u^{*} } \\ {iu_{x} -2\xi u} & {-i\left|u\right|^{2} -2i\xi ^{2} } \end{array}\right).\] 

\begin{lemma}
	Let $f(x,\xi ,t)$ be solution of the equation \eqref{ZEqnNum935108} and let $\phi _{n} (x,\xi ,t),\, \, \, n=1,2,...,2N$ be any functions, which satisfy the conditions
	\begin{equation}\label{eqlem1}
		\frac{\partial \phi _{n} }{\partial x} =G_{n}^{T} f,\, \, \, \, \, n=1,\; 2,\; ...,\; 2N.  
	\end{equation}
	Then, the function $G_{n} (x,t)$ satisfy the equalities
	\begin{equation} \label{ZEqnNum867917} 
		G_{n}^{T} \sigma _{3} f+i(\xi -\xi _{n} )\phi _{n} =0\, ,\, \, \, n=1,\; 2,\; ...,\; 2N 
	\end{equation} 
	and the function 
	\begin{equation} \label{ZEqnNum416287} 
		H=\frac{\partial f}{\partial t} -Af+\sum _{n=1}^{2N}F_{n} \phi _{n}   
	\end{equation} 
	satisfies the equation \eqref{ZEqnNum935108} for any $\xi \in \Sigma $.
\end{lemma}

\subsection{Evolution equation for the scattering data in the case of a source satisfying the conditions (A)}

Let us take matrix Jost solutions $f^{-} (x,\xi ,t)$ and $f^{+} (x,\xi ,t)$ for $\xi \in \Sigma $ as the solution $f(x,\xi ,t)$ and $\xi =\xi _{n} ,\, \, n=1,2,...,N$ are eigenvalues of the equation \eqref{ZEqnNum935108}. According to the definition of eigenfunctions, there are $\alpha _{n}^{} (t)$ and $\beta _{n}^{} (t)$ such that the relations hold
\begin{equation} \label{ZEqnNum492510} 
	F_{n} (x,t)=\alpha _{n} (t)\psi _{n} (x,t),\, \, \, \, G_{n} (x,t)=\beta _{n} (t)\sigma _{1} \varphi _{n} (x,t),\, \, n=1,2,...,N 
\end{equation} 
According to these relations, due to the assumptions \eqref{ZEqnNum767561}, we obtain
\begin{equation} \label{ZEqnNum101214} 
	F_{N+n} (x,t)=\alpha _{n}^{*} (t)\bar{\psi }_{n} (x,t),\, \, G_{N+n} (x,t)=\beta _{n}^{*} (t)\sigma _{1} \bar{\varphi }_{n} (x,t),\, \, \, n=1,2,...,N.  
\end{equation} 

By definition functions $G_{n} (x,t)$, belong to the $L^{2} (R)$ for all $t\ge 0$ and matrix Jost solutions $f^{-} (x,\xi ,t)$, $f^{+} (x,\xi ,t)$ are bounded for all $\xi \in \Sigma $. Therefore $\phi _{n}^{-} \in L^{2} (R)$ and $\phi _{n}^{+} \in L^{2} (R)$ for all $t\ge 0$ and $\xi \in \Sigma $. Hence, by virtue of \eqref{ZEqnNum867917} it follows that at any $\xi \in \Sigma $ and $n=1,2,...,2N$ the asymptotics 
\begin{equation} \label{ZEqnNum380423} 
	\begin{aligned} &\phi _{n}^{-} (x,\xi ,t)\to 0\, \, {\rm as}\, \, x\to -\infty , \\ 
		&\phi _{n}^{+} (x,\xi ,t)\to 0\, \, {\rm as}\, \, x\to \infty  \end{aligned} 
\end{equation} 
are valid. So, from \eqref{eqlem1} for $n=1,2,...,2N$ we obtain the following expressions
\begin{equation} \label{ZEqnNum779844} 
	\phi _{n}^{-} =\int _{-\infty }^{x}G_{n}^{T} (s{\rm ,}t)f^{-} (s,\xi ,t)ds ,\, \, \phi _{n}^{+} =-\int _{x}^{\infty }G_{n}^{T} (s,t)f^{+} (s,\xi ,t)ds .  
\end{equation} 

Using the matrix Jost solutions $f^{+} $ and $f^{-} $ of equation \eqref{ZEqnNum935108}, we rewrite equality  \eqref{ZEqnNum416287} in the form 
\begin{equation} \label{ZEqnNum200928} 
	H^{-} =\frac{\partial f^{-} }{\partial t} -Af^{-} +\sum _{n=1}^{2N}F_{n} \phi _{n}^{-}   
\end{equation} 
and
\begin{equation} \label{ZEqnNum345127} 
	H^{+} =\frac{\partial f^{+} }{\partial t} -Af^{+} +\sum _{n=1}^{2N}F_{n} \phi _{n}^{+}  .  
\end{equation} 
These functions satisfy the equation \eqref{ZEqnNum935108} for any $\xi \in \Sigma $. Therefore, $H^{+} $ and $H^{-} $ are linearly dependent on $f^{+} $ and $f^{-} $, respectively, i.e., there exist such $C_{0}^{-} (\xi ,t)$ and $C_{0}^{+} (\xi ,t)$ that the following identities hold
\[
H^{-} (x,\xi ,t)=f^{-} (x,\xi ,t)C_{0}^{-} (\xi ,t),\, \, \, H^{+} (x,\xi ,t)=f^{+} (x,\xi ,t)C_{0}^{+} (\xi ,t).  
\]
By virtue of the definition of the matrix $A$, from relations \eqref{ZEqnNum200928}, \eqref{ZEqnNum345127} and from asymptotics \eqref{ZEqnNum770859}, \eqref{ZEqnNum380423} we obtain
\begin{align}
	&{H^{-} (x,\xi ,t)\to -(i\rho ^{2} +2i\xi p)E^{-} (x,\xi ,t)\sigma _{3} ,\, \, \, x\to -\infty ,} \\ 
	&{H^{+} (x,\xi ,t)\to -(i\rho ^{2} +2i\xi p)E^{+} (x,\xi ,t)\sigma _{3} ,\, \, \, x\to \infty .} \end{align}
By the uniqueness of the Jost solutions we get
\begin{equation}
	\begin{aligned}\label{ZEqnNum408091}		
		&H^{-} (x,\xi ,t)=-(i\rho ^{2} +2i\xi p)f^{-} (x,\xi ,t)\sigma _{3} , \\ &H^{+} (x,\xi ,t)=-(i\rho ^{2} +2i\xi p)f^{+} (x,\xi ,t)\sigma _{3} .
	\end{aligned} 	
\end{equation}

We introduce the function $\mathcal{H} $ in the following form
\[\mathcal{H} =H_{}^{-} (x,\xi ,t)-H_{}^{+} (x,\xi ,t)S(\xi ,t).\] 
Based on equalities \eqref{ZEqnNum733217} and \eqref{ZEqnNum408091}, the function $\mathcal{H} $ can be rewritten in the form
\begin{equation} \label{ZEqnNum911293} 
	\mathcal{H} =(i\rho ^{2} +2i\xi p)f^{+} (x,\xi ,t)[\sigma _{3} ,S(\xi ,t)].   
\end{equation} 
On the other hand, by virtue of \eqref{ZEqnNum733217}, \eqref{ZEqnNum200928} and \eqref{ZEqnNum345127} the equality 
\begin{align} 
	\mathcal{H} =&H^{-} (x,\xi ,t)-H^{+} (x,\xi ,t)S(\xi ,t)=f^{+} (x,\xi ,t)S_{t} (x,\xi ,t)+ \nonumber\\ 
	+&\sum _{n=1}^{2N}[F_{n} (x,t)\phi _{n}^{-} (x,\xi ,t)-F_{n} (x,t)\phi _{n}^{+} (x,\xi ,t)S(\xi ,t)] 	\label{ZEqnNum777118}
\end{align} 
holds.

Based on equality \eqref{ZEqnNum867917} this relation becomes
\[\mathcal{H} =f^{+} (x,\xi ,t)S_{t} (\xi ,t)+\] 
\[+\sum _{n=1}^{2N}\frac{i}{\xi -\xi _{n} } [F_{n} (x,t)G_{n}^{T} (x,t)\sigma _{3} f^{-} (x,\xi ,t)-F_{n} (x,t)G_{n}^{T} (x,t)\sigma _{3} f^{+} (x,\xi ,t)S(\xi ,t)] .\] 
Finally, based on \eqref{ZEqnNum733217}, we obtain
\begin{equation} \label{ZEqnNum724214} 
	\mathcal{H} =f^{+} (x,\xi ,t)S_{t} (\xi ,t).  
\end{equation} 
Comparing equalities \eqref{ZEqnNum911293} and \eqref{ZEqnNum724214} we have
\[(i\rho ^{2} +2i\xi p)f^{+} (x,\xi ,t)[\sigma _{3} ,S(\xi ,t)]=f^{+} (x,\xi ,t)S_{t} (\xi ,t).\] 
Therefore, for all $\xi \in \Sigma $ we have the relation
\[S_{t} (\xi ,t)-(i\rho ^{2} +2i\xi p)[\sigma _{3} ,S(\xi ,t)]=0,\] 
i.e.
\[\frac{d}{dt} a(\xi ,t)=0,\,\,\,\frac{d}{dt} b(\xi ,t)=-2(i\rho ^{2} +2i\xi p)b(\xi ,t).\] 
Since, the function $a(\xi ,t)$ does not depend on $t$, hence we conclude that its zeros $\xi _{n}$ also do not depend on $t$.

Based on identities \eqref{ZEqnNum200928} and \eqref{ZEqnNum345127}, we write the following equalities
\begin{gather} 
	H_{1}^{-} (x,\xi _{n} ,t)=\frac{\partial \varphi _{m} (x,t)}{\partial t} -A(x,\xi _{n} ,t)\varphi _{m} (x,t)+\sum _{n=1}^{2N}F_{n} (x,t)\phi _{1,n}^{-} (x,\xi _{n} ,t)  \label{ZEqnNum440572} 
	\\
	H_{2}^{+} (x,\xi _{n} ,t)=\frac{\partial \psi _{m} (x,t)}{\partial t} -A(x,\xi _{n} ,t)\psi _{m} (x,t)+\sum _{n=1}^{2N}F_{n} (x,t)\phi _{2,n}^{+} (x,\xi _{n} ,t)  \label{ZEqnNum179643} 
\end{gather} 
By virtue of the definition of the matrix $A$, from relations \eqref{ZEqnNum440572}, \eqref{ZEqnNum179643} and from asymptotics \eqref{ZEqnNum260418}, \eqref{ZEqnNum830052} we obtain
\begin{equation} \label{ZEqnNum863995} 
	\begin{aligned}
		&H_{1}^{-} (x,\xi _{m} ,t)=(-i\rho ^{2} -2i\xi _{m} p_{m} )\varphi _{m} (x,t), \\ &H_{2}^{+} (x,\xi _{m} ,t)=(i\rho ^{2} +2i\xi _{m} p_{m} )\psi _{m} (x,t).  
	\end{aligned}
\end{equation} 
We now introduce the following functions
\[\mathcal{H} _{m} =H_{1}^{-} (x,\xi _{m} ,t)-c_{m} (t)H_{2}^{+} (x,\xi _{m} ,t),\, \, \, m=1,2,...,2N.\] 
Using equalities \eqref{ZEqnNum197225}, \eqref{ZEqnNum863995} the function $\mathcal{H} _{m} $ can be rewritten in the form
\begin{equation} \label{ZEqnNum654365} 
	\mathcal{H} _{m} =(-2i\rho ^{2} -4i\xi _{m} p_{m} )\varphi _{m} (x,t).  
\end{equation} 
Substituting instead of $\phi _{1,n}^{-} (x,\xi ,t)$ and $\phi _{2,n}^{+} (x,\xi ,t)$ the expressions from \eqref{ZEqnNum779844} into equalities \eqref{ZEqnNum440572}, \eqref{ZEqnNum179643} and using \eqref{ZEqnNum197225}, we obtain
\begin{equation} \label{ZEqnNum737308} 
	H_{1}^{-} (x,\xi _{m} ,t)-c_{m} (t)H_{2}^{+} (x,\xi _{m} ,t)=\frac{dc_{m} (t)}{dt} \psi _{m} (x,t)+\sum _{n=1}^{2N}F_{n} (x,t)\int _{-\infty }^{\infty }G_{n}^{T} (s,t) \varphi _{m} (s,t)ds  
\end{equation} 
If $\xi _{m} \ne \xi _{n} $, according to equation \eqref{ZEqnNum867917} we get
\[\int _{-\infty }^{\infty }G_{n}^{T} (s,t) \varphi _{m} (s,t)ds=0.\] 
According to \eqref{ZEqnNum492510} and \eqref{ZEqnNum101214}, equality \eqref{ZEqnNum737308} can be rewritten in the form
\begin{equation} \label{ZEqnNum487278} 
	\begin{gathered} {H_{1}^{-} (x,\xi _{m} ,t)-c_{m} (t)H_{2}^{+} (x,\xi _{m} ,t)=\frac{dc_{m} (t)}{dt} \psi _{m} (x,t)+} \\ {+\left(\int _{-\infty }^{\infty }G_{m}^{T} (s,t)F_{m} (s,t) ds+\int _{-\infty }^{\infty }G_{N+m}^{T} (s,t)F_{N+m} (s,t) ds\right)\varphi _{m} (x,t).} \end{gathered} 
\end{equation} 
Comparing equalities \eqref{ZEqnNum654365} and \eqref{ZEqnNum487278} we obtain
\[(-2i\rho ^{2} -4i\xi _{m} p_{m} )\varphi _{m} (x,t)=\] 
\[=\frac{dc_{m} (t)}{dt} \psi _{m} (x,t)+\left(\int _{-\infty }^{\infty }G_{m}^{T} (s,t)F_{m} (s,t)ds+\int _{-\infty }^{\infty }G_{N+m}^{T} (s,t)F_{N+m} (s,t)ds  \right)\varphi _{m} (x,t).\] 
Finally, using these equalities and taking into account \eqref{eq17} and \eqref{ZEqnNum197225} we determine
\[\frac{dc_{m} (t)}{dt} =(-2i\rho ^{2} -4i\xi _{m} p_{m} -A_{m} (t)-A_{m}^{*} (t))c_{m} (t).\] 

Thus, we have proved the following theorem.

\begin{theorem}\label{thm1}
	If functions $u(x,t),$ $F_{k} (x,t)$, $G_{k} (x,t)$, $k=1,\; 2,\; ...,N$ are the solutions of the problem \eqref{eq11}-\eqref{eq16} in\textbf{ }the case of a source satisfying the conditions\textbf{ }(A), then the scattering data for the equation \eqref{ZEqnNum935108} satisfy the following relations
	\[a(\xi ,t)=a(\xi ,0), \] 
	\[b(\xi ,t)=b(\xi ,0)\exp (-2i\rho ^{2} t-4i\xi pt)\quad  \text{for} \quad  \xi \in \Sigma , \]
	\[\xi _{k} (t)=\xi _{k} (0), \quad k=1,\; 2,\; ...,N.\] 
	\[c_{k} (t)=c_{k} (0)\exp (-2i\rho ^{2} t-4i\xi _{k} p_{k} t-\int _{0}^{t}(A_{k} (\tau )+A_{k}^{*} (\tau ))d\tau  ), \quad k=1,\; 2,\; ...,N.\]
\end{theorem}

\subsection{Evolution equation for the scattering data in the case of a source satisfying the conditions (B)}

Let us take matrix Jost solutions $f^{-} (x,\xi ,t)$ and $f^{+} (x,\xi ,t)$ for $\xi \in \Sigma $ as the solution $f(x,\xi ,t)$ and $\xi =\xi _{n} ,\, \, n=1,2,...,N$ are eigenvalues of the equation \eqref{ZEqnNum935108}. According to the definition of eigenfunctions of equation \eqref{ZEqnNum935108}, there are $\alpha _{n}^{} (t)$ such that the relations
\begin{equation} \label{ZEqnNum594420} 
	F_{n} (x,t)=\alpha _{n} (t)\psi _{n} (x,t),\, \, \, \, \, F_{N+n} (x,t)=\alpha _{n}^{*} (t)\bar{\psi }_{n} (x,t),\, \, \, \, n=1,2,...,N .
\end{equation} 
Due to the assumptions (B) the functions $G_{n} (x,t)$ are unbounded functions. So, there are $\beta _{n} (t)$ such that which follow the equalities
\begin{equation} \label{ZEqnNum360250} 
	\begin{aligned}
		&G_{n} (x,t)=\frac{\beta_{n}(t)}{\dot{a}(\xi_{n},t)}\sigma _{1} \varphi _{n} (x,t)+\sigma _{1} h_{n} (x,t),\\
		&G_{N+n} (x,t)=\frac{\beta _{n}^{*} (t)}{\dot{\bar{a}}(\xi_{n},t)}\sigma _{1} \bar{\varphi }_{n} (x,t)+\sigma _{1} \bar{h}_{n} (x,t),\, \, n=1,2,...,N. 
	\end{aligned} 
\end{equation} 

One can easily see from \eqref{eq18} and \eqref{ZEqnNum819680}, that the quantities $\alpha_{n}(t)$  satisfy the following equalities 
\begin{equation} \label{ZEqnNum951420} 
	\alpha _{n} (t) =-\frac{\rho ^{2} }{2p_{n} (\xi _{n} -p_{n} )} B_{n} (t),\, \, \, \, \alpha _{n}^{*}(t)  =\frac{\rho ^{2} }{2p_{n} (\xi _{n} +p_{n} )} B_{n}(t),\, \, n=1,2,...,N.  
\end{equation} 

Using equalities \eqref{ZEqnNum197225}, \eqref{ZEqnNum867917}, \eqref{ZEqnNum594420}, \eqref{ZEqnNum360250} and asymptotics \eqref{ZEqnNum260418}, \eqref{ZEqnNum922580}
we can verify that at any $\xi \in \Sigma $ and when $x\to \infty $ the following asymptotics are valid:
\[
F_{n} \phi _{n}^{+} \sim \frac{i\alpha _{n}(t)}{\xi -\xi _{n} } \left(\begin{array}{cc} {\frac{(\xi _{n} -p_{n} )^{2} }{\rho ^{2} } } & {-\frac{i(\xi _{n} -p_{n} )}{\rho } e^{-i\alpha _{+} +2i\rho ^{2} t} } \\ {\frac{i(\xi _{n} -p_{n} )}{\rho } e^{i\alpha _{+} -2i\rho ^{2} t} } & {1} \end{array}\right)\sigma _{3} E^{+} (x,\xi ,t),  
\]
\[
F_{N+n} \phi _{N+n}^{+} \sim \frac{i\alpha _{n}^{*}(t) }{\xi -\xi _{n} } \left(\begin{array}{cc} {1} & {-\frac{i(\xi _{n} +p_{n} )}{\rho } e^{-i\alpha _{+} +2i\rho ^{2} t} } \\ {\frac{i(\xi _{n} +p_{n} )}{\rho } e^{i\alpha _{+} -2i\rho ^{2} t} } & {\frac{(\xi _{n} +p_{n} )^{2} }{\rho ^{2} } } \end{array}\right)\sigma _{3} E^{+} (x,\xi ,t).  
\]
Taking into account of equalities \eqref{ZEqnNum197225}, \eqref{ZEqnNum867917}, \eqref{ZEqnNum594420}, \eqref{ZEqnNum360250} and asymptotics \eqref{ZEqnNum260418}, \eqref{ZEqnNum922580} we are convinced that at any $\xi \in \Sigma $ and $x\to -\infty $ there hold the asymptotics
\[
F_{n} \phi _{n}^{-} \sim -\frac{i\alpha _{n} (t) }{\xi -\xi _{n} } \left(\begin{array}{cc} {1} & {-\frac{i(\xi _{n} -p_{n} )}{\rho } e^{-i\alpha _{-} +2i\rho ^{2} t} } \\ {\frac{i(\xi _{n} -p_{n} )}{\rho } e^{i\alpha _{-} -2i\rho ^{2} t} } & {\frac{(\xi _{n} -p_{n} )^{2} }{\rho ^{2} } } \end{array}\right)\sigma _{3} E^{-} (x,\xi ,t),  
\] 
\[
F_{N+n} \phi _{N+n}^{-} \sim -\frac{i\alpha _{n}^{*}(t)  }{\xi -\xi _{n} } \left(\begin{array}{cc} {\frac{(\xi _{n} +p_{n} )^{2} }{\rho ^{2} } } & {-\frac{i(\xi _{n} +p_{n} )}{\rho } e^{-i\alpha _{-} +2i\rho ^{2} t} } \\ {\frac{i(\xi _{n} +p_{n} )}{\rho } e^{i\alpha _{-} -2i\rho ^{2} t} } & {1} \end{array}\right)\sigma _{3} E^{-} (x,\xi ,t).  
\]
Using equalities \eqref{ZEqnNum951420} one can easily verify that at any $\xi \in \Sigma $ the asymptotic 
\[
F_{n} \phi _{n}^{+} +F_{N+n} \phi _{N+n}^{+} \to 0 \quad \mathrm{for}  \quad x\to \infty,
\]
and
\[F_{n} \phi _{n}^{-} +F_{N+n} \phi _{N+n}^{-} \to 0 \quad \mathrm{for}   \quad x\to -\infty 
\]
are valid. 

Hence, it follows that the quantities $H^{-} (x,\xi ,t)$ and $H^{+} (x,\xi ,t)$ determined by \eqref{ZEqnNum200928} and \eqref{ZEqnNum345127} satisfy equalities 
\begin{equation} \label{ZEqnNum355034} 
	\begin{gathered} {H^{-} (x,\xi ,t)=f^{-} (x,\xi ,t)(-i\rho ^{2} -2i\xi p)\sigma _{3} ,\, \, \, } \\ {H^{+} (x,\xi ,t)=f^{+} (x,\xi ,t)(-i\rho ^{2} -2i\xi p)\sigma _{3} .} \end{gathered}  
\end{equation} 
Now, consider the function $\mathcal{H}_m $ of the form
\[\mathcal{H}_m =H^{-} (x,\xi ,t)-H^{+} (x,\xi ,t)S(\xi ,t).\] 
Taking into account \eqref{ZEqnNum355034} we find that
\begin{equation} \label{ZEqnNum469933} 
	\mathcal{H}_m =(i\rho ^{2} +2i\xi p)f^{+} (x,\xi ,t)[\sigma _{3} ,S(\xi ,t)].  
\end{equation} 
From equalities \eqref{ZEqnNum200928}, \eqref{ZEqnNum345127} and \eqref{ZEqnNum733217}  it is easy to get that
\begin{align} 
	&H^{-} (x,\xi ,t)-H^{+} (x,\xi ,t)S(\xi ,t)=f^{+} (x,\xi ,t)S_{t} (x,\xi ,t)+\\ \nonumber
	&+\sum _{n=1}^{2N}[F_{n} (x,t)\phi _{n}^{-} (x,\xi ,t)-F_{n} (x,t)\phi _{n}^{+} (x,\xi ,t)S(\xi ,t)] .  \label{GrindEQ__1_69_} 
\end{align} 
Using equalities \eqref{ZEqnNum867917}, we obtain 
\[H^{-} (x,\xi ,t)-H^{+} (x,\xi ,t)S(\xi ,t)=f^{+} (x,\xi ,t)S_{t} (\xi ,t)+\] 
\[+\sum _{n=1}^{2N}\frac{i}{\xi -\xi _{n} } [F_{n} (x,t)G_{n}^{T} (x,t)\sigma _{3} f^{-} (x,\xi ,t)-F_{n} (x,t)G_{n}^{T} (x,t)\sigma _{3} f^{+} (x,\xi ,t)S(\xi ,t)] .\] 
By virtue of  \eqref{ZEqnNum733217}, it follows that
\begin{equation} \label{ZEqnNum157313} 
	H^{-} (x,\xi ,t)-H^{+} (x,\xi ,t)S(\xi ,t)=f^{+} (x,\xi ,t)S_{t} (\xi ,t) 
\end{equation} 
Comparing equalities \eqref{ZEqnNum469933} and \eqref{ZEqnNum157313} we have
\[(i\rho ^{2} +2i\xi p)f^{+} (x,\xi ,t)[\sigma _{3} ,S(\xi ,t)]=f^{+} (x,\xi ,t)S_{t} (\xi ,t).\] 
Therefore, for all $\xi \in \Sigma $ we have 
\[S_{t} (\xi ,t)-(i\rho ^{2} +2i\xi p)[\sigma _{3} ,S(\xi ,t)]=0,\] 
i.e.
\[\frac{d}{dt} a(\xi ,t)=0,\,\,\,\frac{d}{dt} b(\xi ,t)=-2(i\rho ^{2} +2i\xi p)b(\xi ,t).\] 
Thus, we conclude that the function $a(\xi ,t)$ does not depend on $t$, so the zeros $\xi _{n} (t)$ of function $a(\xi ,t)$ do not depend on $t$.

Let us now find the evolution of the normalizing constants $c_{m} (t),\, \, m=1,2,...,N$.

We now introduce the following functions
\begin{equation} \label{ZEqnNum112089} 
	\mathcal{H}_m =H_{1}^{-} (x,\xi _{m} ,t)-c_{m} (t)H_{2}^{+} (x,\xi _{m} ,t),\, \, m=1,2,...,N,  
\end{equation} 
where
\begin{gather} 
	H_{1}^{-} (x,\xi _{m} ,t)=\frac{\partial \varphi _{m} (x,t)}{\partial t} -A(x,\xi _{m} ,t)\varphi _{m} (x,t)+\sum _{n=1}^{2N}\Phi _{n} (x,t)\varphi _{1,n}^{-} (x,\xi _{m} ,t) , \label{ZEqnNum662125}  
	\\ 
	H_{2}^{+} (x,\xi _{m} ,t)=\frac{\partial \psi _{m} (x,t)}{\partial t} -A(x,\xi _{m} ,t)\psi _{m} (x,t)+\sum _{n=1}^{2N}\Phi _{n} (x,t)\varphi _{2,n}^{+} (x,\xi _{m} ,t) .  \label{ZEqnNum300520} 
\end{gather} 
It is easy to show that
\begin{equation} \label{ZEqnNum978528} 
	H_{1}^{-} (x,\xi _{m} ,t)=(-i\rho ^{2} -2i\xi _{m} p_{m} )\varphi _{m} (x,t),\, \, H_{2}^{+} (x,\xi _{m} ,t)=(i\rho ^{2} +2i\xi _{m} p_{m} )\psi _{m} (x,t).  
\end{equation} 
Substituting \eqref{ZEqnNum978528} into \eqref{ZEqnNum112089} and using equalities \eqref{ZEqnNum197225}, we get for $m=1,2,...,N$
\begin{equation} \label{ZEqnNum227254} 
	\mathcal{H}_m =(-2i\rho ^{2} -4i\xi _{m} p_{m} )\varphi _{m} (x,t).  
\end{equation} 
On the other hand, using equalities \eqref{ZEqnNum662125}, \eqref{ZEqnNum300520} and \eqref{ZEqnNum197225} we obtain 
\[H_{1}^{-} (x,\xi _{m} ,t)-c_{m} (t)H_{2}^{+} (x,\xi _{m} ,t)=\frac{dc_{m} (t)}{dt} \psi _{m} (x,t)+\] 
\[+\sum _{{}^{n=1}_{n\ne m\, }}^{2N}\frac{i}{\xi _{m} -\xi _{n} } [F_{n} (x,t)G_{n}^{T} (x,t)\sigma _{3} \varphi _{m} (x,t)-c_{m} (t)F_{n}^{+} (x,t)G_{n}^{T} (x,t)\sigma _{3} \psi _{m} (x,t)] +\] 
\[+iF_{m} (x,t)G_{m}^{T} (x,t)\sigma _{3} \frac{\partial }{\partial \xi } \left. \left(\varphi (x,\xi ,t)-c_{m} (t)\psi (x,\xi ,t)\right)\right|_{\xi =\xi _{m} } +\] 
\[+iF_{N+m} (x,t)G_{N+m}^{T} (x,t)\sigma _{3} \frac{\partial }{\partial \xi } \left. \left(\varphi (x,\xi ,t)-c_{m} (t)\psi (x,\xi ,t)\right)\right|_{\xi =\xi _{m} } .\] 
According to \eqref{ZEqnNum197225} and \eqref{ZEqnNum510514}, this equation can be rewritten in the following form
\[H_{1}^{-} (x,\xi _{m} ,t)-c_{m} (t)H_{2}^{+} (x,\xi _{m} ,t)=\frac{dc_{m} (t)}{dt} \psi _{m} (x,t)+\] 
\[+i\dot{a}(\xi _{n} ,t)F_{m} (x,t)G_{m}^{T} (x,t)\sigma _{3} h_{m} (x,t)+i\dot{a}(\xi _{n} ,t)F_{N+m} (x,t)G_{N+m}^{T} (x,t)\sigma _{3} h_{m} (x,t).\] 
Further, by virtue \eqref{eq18}, \eqref{ZEqnNum594420} and \eqref{ZEqnNum360250}, we obtain the equality
\begin{align} 
	&H_{1}^{-} (x,\xi _{m} ,t)-c_{m} (t)H_{2}^{+} (x,\xi _{m} ,t)=\frac{dc_{m} (t)}{dt} \psi _{m} (x,t)- \nonumber\\ 
	&-i\beta _{m} (t) B_{m} (t)\varphi _{m} (x,t)+i\beta _{m}^{*} (t) B_{m} (t){\varphi }_{m} (x,t). \label{ZEqnNum889349} 
\end{align} 
Comparing equalities \eqref{ZEqnNum227254} and \eqref{ZEqnNum889349}, we obtain
\[(-2i\rho ^{2} -4i\xi _{m} p_{m} )\varphi _{m} (x,t)=\frac{dc_{m} (t)}{dt} \psi _{m} (x,t)-i\beta _{m} (t) B_{m} (t)\varphi _{m} (x,t)+i\beta _{m}^{*} (t) B_{m} (t){\varphi }_{m} (x,t)\] 
hence, taking into account equalities \eqref{ZEqnNum197225},  we find
\[\frac{dc_{m} (t)}{dt} =(-2i\rho ^{2} -4i\xi _{m} p_{m} +i(\beta _{m} (t)  -\beta _{m}^{*} (t)) B_{m} (t) )c_{m} (t).\] 

Thus, we have proved the following theorem.

\begin{theorem}\label{thm2} If functions $u(x,t)$, $F_{k} (x,t)$, $G_{k} (x,t)$, $k=1,\; 2,\; ...,N$ are the  solutions of the problem \eqref{eq11}-\eqref{eq16} in the case of a source satisfying the conditions (B), then the scattering data for the equation \eqref{ZEqnNum935108} satisfy the following relations
	\[a(\xi ,t)=a(\xi ,0), \] 
	\[b(\xi ,t)=b(\xi ,0)\exp (-2i\rho ^{2} t-4i\xi pt) \quad \text{for} \quad \xi \in \Sigma , \]
	\[\xi _{k} (t)=\xi _{k} (0), \quad k=1,\; 2,\; ...,N.\] 
	\[c_{k} (t)=c_{k} (0)\exp (-2i\rho ^{2} -4i\xi _{k} p_{k} +i\int _{0}^{t}\left(\beta _{k} (\tau ) -\beta _{k}^{*} (\tau ) \right) B_{k} (\tau )d\tau , k=1,\; 2,\; ...,N.\] 
\end{theorem}

We will illustrate inverse scattering method of constructing exact solutions to the NLS equation with concrete example.

\begin{example}\label{exm1}
	Let the initial function $u_{0} (x)$ have the form
	\[u_{0} (x)=\rho \cdot \frac{e^{i\alpha _{+} } e^{\nu x} +e^{i\alpha _{-} } e^{-\nu x} }{e^{\nu x} +ce^{-\nu x} } .\] 
	Where $\alpha _{+} ,\, \alpha _{-} ,\, \rho ,\, \nu $, $c$ are positive real numbers and $\rho >\nu $.

\end{example}

In this case, solving the direct scattering problem for the equation \eqref{ZEqnNum935108}, we obtain 
\[a(\xi ,0)=\frac{\xi +p-\zeta -i\nu }{\xi +p-\zeta +i\nu } ,\, \, \zeta =\sqrt{\rho ^{2} -\nu ^{2} } , \] 
\[b(\xi ,0)=0, \quad \xi _{1} (0)=\zeta , \quad c_{1} (0)=\frac{i(\zeta -i\nu )}{\rho } ce^{i\alpha _{-} } .\] 

Based on Theorem \ref{thm1}, we can show the evolution of the scattering data in the following form
\[a(\xi ,t)=\frac{\xi +p-\zeta -i\nu }{\xi +p-\zeta +i\nu } ,\, \, \zeta =\sqrt{\rho ^{2} -\nu ^{2} } , \] 
\[b(\xi ,t)=0, \quad \xi _{1} (t)=\zeta , \] 
\[c_{1} (t)=\frac{i(\zeta -i\nu )}{\rho } \cdot c\cdot \exp (i\alpha _{-} -2i\rho ^{2} t+4\zeta \nu t-\int _{0}^{t}(A_{k} (\tau )+A_{k}^{*} (\tau ))d\tau  ).\] 

Applying the procedure of the inverse scattering problem, we find
\[u(x,t)=\rho e^{-2i\rho ^{2} t} \cdot \frac{e^{i\alpha _{+} } e^{\nu x} +e^{i\alpha _{-} } ce^{-\nu x+4\zeta \nu t-g(t)} }{e^{\nu x} +ce^{-\nu x+4\zeta \nu t-g(t)} } ,\] 
where \textbf{$g(t)=\int _{0}^{t}(A_{1} (\tau )+A_{1}^{*} (\tau )) \, d\tau $}.

Using representation \eqref{ZEqnNum820131} and conditions \eqref{eq17}, we obtain
\[ {F_{1} =\alpha _{1}(t) \cdot \left(\begin{array}{c} {-\frac{i(\zeta -i\nu )}{\rho } \cdot e^{-i\alpha _{+} +2i\rho ^{2} t} } \\ {1} \end{array}\right)\cdot \frac{1}{e^{\nu x} +ce^{-\nu x+4\zeta \nu t-g(t)} } ,} \]
\[ {G_{1} =\frac{\nu A_{1} (t)}{\alpha _{1}(t) } \cdot \left(\begin{array}{c} {\frac{i(\zeta -i\nu )}{\rho } \cdot e^{i\alpha _{-} -2i\rho ^{2} t} } \\ {1} \end{array}\right)\cdot \frac{ce^{4\zeta \nu t-g(t)} }{e^{\nu x} +ce^{-\nu x+4\zeta \nu t-g(t)} } .} \] 
Analogously, in the case (B), using results of Theorem \ref{thm2}, we obtain 
\[u(x,t)=\rho e^{-2i\rho ^{2} t} \cdot \frac{e^{i\alpha _{+} } e^{\nu x} +e^{i\alpha _{-} } ce^{-\nu x+4\zeta \nu t-g(t)} }{e^{\nu x} +ce^{-\nu x+4\zeta \nu t-g(t)} } ,\] 
\[F_{1} =-\frac{\rho^2 B_1(t)}{2i\nu(\zeta-i\nu)} \cdot \left(\begin{array}{c} {-\frac{i(\zeta -i\nu )}{\rho } \cdot e^{-i\alpha _{+} +2i\rho ^{2} t} } \\ {1} \end{array}\right)\cdot \frac{1}{e^{\nu x} +ce^{-\nu x+4\zeta \nu t-g(t)} } ,\] 
\[ {G_{1} =\frac{-2\nu }{\zeta+i\nu}\cdot \left(\nu\beta_{1}(t)-2x\zeta+i\sigma_3\right) \cdot \left(\begin{array}{c} {\frac{i(\zeta -i\nu )}{\rho } \cdot e^{i\alpha _{-} -2i\rho ^{2} t} } \\ {1} \end{array}\right)\cdot \frac{ce^{4\zeta \nu t-g(t)} }{e^{\nu x} +ce^{-\nu x+4\zeta \nu t-g(t)} } +} \] 
\[
+\left(\begin{array}{c} {\frac{i(\zeta -i\nu )}{\rho } \cdot e^{i\alpha _{+} -2i\rho ^{2} t} } \\ {1} \end{array}\right)\cdot \frac{e^{2\nu x} }{e^{\nu x} +ce^{-\nu x+4\zeta \nu t-g(t)} } + 
\]
\[ 
{+\left(\begin{array}{c} {-\frac{i(\zeta -i\nu )}{\rho } \cdot e^{i\alpha _{-} -2i\rho ^{2} t}} \\ {e^{i\theta} } \end{array}\right)\cdot \frac{c^{2} e^{-2\nu x+8\zeta \nu t-2g(t) } }{e^{\nu x} +ce^{-\nu x+4\zeta \nu t-g(t)} } }, 
\] 
where $g(t)=-i\int _{0}^{t}\left(\beta _{1} (\tau ) -\beta _{1}^{*} (\tau ) \right) B_{1} (\tau )d\tau $.

%%%% Acknowledgments %%%%%%%%
\section*{Acknowledgments}
This research was supported by program ``Short-term research internships of young scientists in leading foreign scientific organizations'' of the Ministry of Innovative Development of the Republic of Uzbekistan. Endless gratitude to Professor Rogrigo Lopez for the great support of my research visit to University of Santiago de Compostela.

%%%% Bibliography  %%%%%%%%%%


\begin{thebibliography}{10}
	\footnotesize\itemsep=0pt
	\providecommand{\url}[1]{#1}
	\providecommand{\urlprefix}{}
	
	
	\bibitem{Ablowitz1974}
	Ablowitz M.J., Kaup D.J., Newell A.C., Segur H., The inverse scattering
	transform-{F}ourier analysis for nonlinear problems, \textit{Studies in Appl.
		Math.} \textbf{53} (1974), 249--315.
	
	\bibitem{Biondini2016}
	Biondini G., Fagerstrom E., Prinari B., Inverse scattering transform for the
	defocusing nonlinear {S}chr\"{o}dinger equation with fully asymmetric
	non-zero boundary conditions, \textit{Phys. D} \textbf{333} (2016), 117--136.
	
	\bibitem{Biondini2014}
	Biondini G., Kova\v{c}i\v{c} G., Inverse scattering transform for the focusing
	nonlinear {S}chr\"{o}dinger equation with nonzero boundary conditions,
	\textit{J. Math. Phys.} \textbf{55} (2014), 031506, 22.
	
	\bibitem{Demontis2013}
	Demontis F., Prinari B., van~der Mee C., Vitale F., The inverse scattering
	transform for the defocusing nonlinear {S}chr\"{o}dinger equations with
	nonzero boundary conditions, \textit{Stud. Appl. Math.} \textbf{131} (2013),
	1--40.
	
	\bibitem{Demontis2014}
	Demontis F., Prinari B., van~der Mee C., Vitale F., The inverse scattering
	transform for the focusing nonlinear {S}chr\"{o}dinger equation with
	asymmetric boundary conditions, \textit{J. Math. Phys.} \textbf{55} (2014),
	101505, 40.
	
	\bibitem{Demontis2014Acta}
	Demontis F., van~der Mee C., Characterization of scattering data for the {AKNS}
	system, \textit{Acta Appl. Math.} \textbf{131} (2014), 29--47.
	
	\bibitem{Demontis2022}
	Demontis F., van~der Mee C., A matrix {S}chr\"{o}dinger approach to focusing
	nonlinear {S}chr\"{o}dinger equations with nonvanishing boundary conditions,
	\textit{J. Nonlinear Sci.} \textbf{32} (2022), Paper No. 57, 29.
	
	\bibitem{Faddeev1987}
	Faddeev L.D., Takhtajan L.A., Hamiltonian methods in the theory of solitons,
	Springer Series in Soviet Mathematics, Springer-Verlag, Berlin, 1987.
	
	\bibitem{Khasanov2009}
	Khasanov A.B., Re\u{\i}imberganov A.A., On the finite-density solution of the
	nonlinear {S}chr\"{o}dinger equation with a self-consistent source,
	\textit{Uzbek. Mat. Zh.}  (2009), 123--130.
	
	\bibitem{Melnikov1991}
	Melnikov V.K., Integration of the nonlinear {S}chroedinger equation with a
	self-consistent source, \textit{Comm. Math. Phys.} \textbf{137} (1991),
	359--381.
	
	\bibitem{Melnikov1992}
	Melnikov V.K., Integration of the nonlinear {S}chr\"{o}dinger equation with a
	source, \textit{Inverse Problems} \textbf{8} (1992), 133--147.
	
	\bibitem{Prinari2018}
	Prinari B., Demontis F., Li S., Horikis T.P., Inverse scattering transform and
	soliton solutions for square matrix nonlinear {S}chr\"{o}dinger equations
	with non-zero boundary conditions, \textit{Phys. D} \textbf{368} (2018),
	22--49.
	
	\bibitem{Reyimberganov2021}
	Reyimberganov A.A., Rakhimov I.D., The soliton solutions for the nonlinear
	{S}chr\"{o}dinger equation with self-consistent source, \textit{Izv. Irkutsk.
		Gos. Univ. Ser. Mat.} \textbf{36} (2021), 84--94.
	
	\bibitem{Urazboev2021}
	Urazboev G.U., Reyimberganov A.A., Babadjanova A.K., Integration of the matrix
	nonlinear {S}chr\"{o}dinger equation with a source, \textit{Izv. Irkutsk.
		Gos. Univ. Ser. Mat.} \textbf{37} (2021), 63--76.
	
	\bibitem{Yakhshimuratov2011}
	Yakhshimuratov A., The nonlinear {S}chr\"{o}dinger equation with a
	self-consistent source in the class of periodic functions, \textit{Math.
		Phys. Anal. Geom.} \textbf{14} (2011), 153--169.
	
	\bibitem{Zakharov1971}
	Zakharov V.E., Shabat A.B., Exact theory of two-dimensional self-focusing and
	one-dimensional self-modulation of waves in nonlinear media, \textit{Soviet
		Journal of Experimental and Theoretical Physics} \textbf{61} (1971),
	118--134.
	
	\bibitem{Zakharov1972}
	Zakharov V.E., Shabat A.B., Interaction between solitons in a stable medium,
	\textit{Soviet Journal of Experimental and Theoretical Physics} \textbf{37}
	(1973), 823--828.
	
\end{thebibliography}
\end{document}